\documentclass[10pt, leqno]{article}

\usepackage{mathrsfs}
\usepackage{amsmath}
\usepackage{amssymb}
\usepackage{amsthm}
\usepackage{comment}
\usepackage[dvipdfmx]{hyperref}

\usepackage{accents}

\newtheorem{them}{Theorem}[section]
\newtheorem{prop}{Proposition}[section]

\newtheorem{cor}{Corollary}[section]

\numberwithin{equation}{section}

\def\ep{\epsilon}
\def\dif{\partial}

\begin{document}

\title{
A question on the Cauchy problem in the Gevrey classes for weakly hyperbolic equations}

\author{Tatsuo Nishitani\footnote{Department of Mathematics, Osaka University:  
nishitani@math.sci.osaka-u.ac.jp
}}

\date{}
\maketitle

\def\dif{\partial}
\def\al{\alpha}
\def\be{\beta}
\def\ga{\gamma}
\def\om{\omega}
\def\lam{\lambda}
\def\tika{{\tilde \nu}}
\def\baka{{\bar \nu}}
\def\varep{\varepsilon}
\def\tal{{\tilde\alpha}}
\def\tbe{{\tilde\beta}}
\def\tis{{\tilde s}}
\def\bas{{\bar s}}
\def\R{{\mathbb R}}
\def\N{{\mathbb N}}
\def\C{{\mathbb C}}
\def\Q{{\mathbb Q}}
\def\Ga{\Gamma}
\def\La{\Lambda}
\def\lr#1{\langle{#1}\rangle}
\def\lrD{\langle{D}\rangle}
\def\mD{\lr{ D}_{\mu}}
\def\xim{\lr{\xi}_{\mu}}
\def\co{{\mathcal C}}
\def\op#1#2{{\rm op}^{#1}({#2})}


\begin{abstract}
For a homogeneous polynomial $p$ in $\xi\in\R^n$ with Gevrey coefficients, it is known that the Cauchy problem for any realization of $p$ is well-posed in the  Gevrey class of order $s<2$ if the characteristic roots are real.  In this note, we give examples showing the situation of the converse direction, in particular the optimality of the Gevrey order $s=2$. 
\end{abstract}

\section{Introduction}

Consider a polynomial in $\xi=(\xi_1,\ldots, \xi_n)$ of degree $m$ with Gevrey coeffcietnts
\[
P(x,\xi)=\sum_{|\al|\leq m }a_{\al}(x)\xi^m=p(x, \xi)+\sum_{j=0}^{m-1}P_j(x, \xi)
\]
where $a_{\al}(x)$ are in some Gevrey classes $\gamma^{(s)}(\Omega)$ or $\gamma^{\lr{s}}(\Omega)$ defined in a neighborhood of the origin of $\R^n$ and $p(x,\xi)$, $P_j(x, \xi)$   denotes the homogeneous part of degree $m$ and $j<m$ in $\xi$ respectively.  Gevrey class $\gamma^{(s)}(\Omega)$ is the set of functions $f\in C^{\infty}(\Omega)$ such that for any compact set $K\Subset \Omega$ there exist constants $C>0, A>0$ for which the following inequalities hold:
\begin{equation}
\label{eq:G:teigi}
|D^{\al}f(x)|\leq CA^{|\al|}(|\al|!)^s,\;\; x\in K,\;\;\al\in \N^{n}.
\end{equation}
%
For a given $P(x,\xi)$ we consider its several realizations (quantizations)  as differential operators. We define $\op{t}{P}$, $0\leq t\leq 1$ by
\[
(\op{t}{P}u)(x)=(2\pi)^{-n}\int e^{i(x-y)\xi}P((1-t)x+ty,\xi)u(y)dyd\xi.
\]
Note that, assuming that $a_{\al}(x)$ are constant outside some compact neighborhood of the origin for simplicity, we see
\[
\sum_{|\al|\leq m}a_{\al}(x)D^{\al}u(x)=\op{0}{P}u(x),\quad
\sum_{|\al|\leq m}D^{\al}\big(a_{\al}(x)u(x)\big)=\op{1}{P}u(x).
\]
When $t=1/2$ the quantization $\op{1/2}{P}$ is called Weyl quantization and also denoted by $\op{w}{P}$. Note that
\[
\lim_{\lambda\to\infty}\lambda^{-m}e^{-\lambda x \xi}\op{t}{P}e^{\lambda x \xi}=p(x, \xi)
\]
so that the {\it principal symbol} of $\op{t}{P}$ is independent of realization and given by the highest homogeneous part of $P(x, \xi)$.
 
Denote $x=(x_1, x_2,\ldots, x_n)=(x_1, x')$ and 
consider the Cauchy problem 
\begin{equation}
\label{eq:CPm}
\left\{\begin{array}{ll}
\op{t}{P} u(x)=0,\quad (x_1, x')\in \omega\cap\{x_1>\tau\},\\[8pt]
D_1^ju(0,x')=u_j(x'),\;\; j=0,\ldots, m-1,\;\;x'\in \omega\cap\{x_1=\tau\}
\end{array}\right.
\end{equation}
where $\omega\subset\Omega$ is some open neighborhood of the origin of $\R^{n}$. The Cauchy problem \eqref{eq:CPm} is (uniformly) well-posed in $\gamma^{(s)}$ 
near the origin if there exist $\omega$ and $\ep>0$ such that for any $u_j\in \gamma^{(s)}(\R^{n-1})$ 
and for any $|\tau|<\ep$ the Cauchy problem \eqref{eq:CPm} has a unique solution $u\in C^m(\omega)$. We say that the Cauchy problem is locally solvable in $\gamma^{(s)}$ at the origin if for any $u_j\in \gamma^{(s)}(\R^{n-1})$ one can find a neighborhood $\omega$ of the origin, which may depend on $\{u_j\}$,  such that \eqref{eq:CPm} with $\tau=0$ has a solution $u\in C^m(\omega)$. We assume that the hyperplanes $x_1=\text{const.,}$ are non-characteristic, that is
\begin{equation}
\label{eq:nonchar}
p(x, \theta)\neq 0,\;\; \theta=(1, 0, \ldots, 0),\;\;x\in \Omega,
\end{equation}
which is almost necessary for $C^{\infty}$ well-posedness of the Cauchy problem (\cite{Mi}). Then without restrictions one may assume $a_{(m, 0,\ldots, 0)}(x)=1$. If the Cauchy problem \eqref{eq:CPm} is locally solvable in $\gamma^{(s)}$, $s>1$ at the origin then $p(0, \xi_1,\xi')=0$ has only real roots for any $\xi'\in \R^{n-1}$ (\cite{Ni:1}) so we assume without restrictions 
\begin{equation}
\label{eq:hyperb}
p(x, \xi-i\theta)= 0,\;\;\xi\in \R^{n-1},\;\;x\in\Omega.
\end{equation}
%
The next results are implicit in \cite{Br} and \cite{Ka}.

\begin{them}
\label{thm:hiseiji}Assume $P(x,\xi)=p(x,\xi)+\sum_{j=0}^{r}P_j(x, \xi)$ and the coefficients $a_{\al}(x)$ belong to $\gamma^{(m/r)}$. 
Then for any $0\leq t\leq 1$ the Cauchy problem for $\op{t}{P}$ is well-posed in $\gamma^{(s)}$ near the origin for $1<s<\min{\{m/r, 2\}}$. 
\end{them}
\begin{cor}
\label{cor:seiji}
Assume that the coefficients $a_{\al}(x)$ belong to $\gamma^{(2)}$. 
Then for any $0\leq t\leq 1$ the Cauchy problem for $\op{t}{p}$ is well-posed in $\gamma^{(s)}$ near the origin  for $1<s<2$. 
\end{cor}
To confirm the results, note that under the assumption there is some $M>0$ such that
\begin{equation}
\label{eq:Larsson}
P(x, \xi+i\tau\theta)\neq 0,\quad |\tau|\geq M(1+|\xi|)^{r/m},\quad x,\;\xi\in\R^n
\end{equation}
that is $P$ is $m/r$ -\,hyperbolic (see \cite{Lar}). Then Theorem \ref{thm:hiseiji} was proved for $\op{0}{P}$ in \cite{Ka} and Corollary \ref{cor:seiji} for $\op{0}{p}$ is implicit in \cite{Br}. Next, we recall a formula for change of quantization (e.g. \cite{Mar}). One can pass from any $t$-\,quantization to the $t'$-\,quantization by
\begin{equation}
\label{eq:henko}
\op{t'}{a_{t'}}=\op{t}{a_t},\quad a_{t'}(x,\xi)=e^{-i(t'-t)D_{x}D_{\xi}}a_t(x,\xi)
\end{equation}
for $a(x, \xi)\in S^m_{1, 0}$, symbols of classical pseudodifferential operators. In particular, we have
\[
\op{t}{P(x,\xi)}={\rm op}^0\big(e^{itD_xD_{\xi}}P(x, \xi)\big).
\]
On the other hand, from \cite[Proposition 3]{Br} one has
\[
\big|(\dif_x^{\be}\dif_{\xi}^{\al}p)(x,\xi-i\theta)\big|\leq C_{\al\be}(1+|\xi|)^{|\be|}|p(x,\xi-i\theta)|,\quad \al,\be\in \N^n
\]
which is sufficient to estimate new terms that appear by operating  $e^{itD_xD_{\xi}}$ to $P(x,\xi)$.

\section{A question on Theorem \ref{thm:hiseiji}}

If $P$ is of constant coefficients, $P(x,\xi)=P(\xi)$, the Cauchy problem for $P(D)$ is $\gamma^{(s)}$ well-posed for $1<s<m/r$ (\cite{Lar}) if \eqref{eq:Larsson} holds  where it is understood that the Cauchy problem is $C^{\infty}$ well-posed when $r=0$, which corresponds to the hyperbolicity in the sense of G\aa rding (\cite{Gar}). It is  clear that the results are optimal considering examples $P(D)=D_1^m+cD_n^r$ with a suitable $c\in\C$.  In the variable coefficient case, on the other hand,  we are restricted to $1<s<\min{\{m/r, 2\}}$ in both Theorem \ref{thm:hiseiji} and Corollary \ref{cor:seiji}. Here we give an example showing that one can not exceed $2$, at least when $m\geq 3$.  Consider
\begin{equation}
\label{eq:P:mod}
P_{b}(x, \xi)=\xi_1^3-(\xi_2^2+x_2^2 \xi_n^2)\xi_1-b\, x_2^3\xi_n^3
\end{equation}
with $b\in\R$ which was studied in \cite{BeNi:2}. Note that \eqref{eq:hyperb} is equivalent to $b^2\leq 4/27$. In \cite{BeNi:2} it was proved that there is $0<{\bar b}<2/3\sqrt{3}$ such that the Cauchy problem for $\op{0}{P_{\,{\bar b}}}$ is not locally solvable at the origin in $\gamma^{(s)}$ for $s>2$. 
\begin{prop}
\label{pro:m:3} Let $m\geq 3$ and $n\geq 3$ and consider
\[
p(x, \xi)=\xi_1^{m-3}\big(\xi_1^3-(\xi_2^2+x_2^2\xi_n^2)\xi_1-{\bar b}\,x_2^3\xi_n^3\big)
\]
which is a homogeneous polynomial in $\xi$ of degree $m$ with polynomial coefficients. For any $0\leq t\leq 1$ the Cauchy problem for $\op{t}{p}$ is not locally solvable at the origin in $\gamma^{(s)}$ for $s>2$, in particular, ill-posed in $\gamma^{(s)}$, $s>2$ near the origin.
\end{prop}
From \eqref{eq:henko} one sees that $\op{t'}{p}=\op{t}{p}$ for any $0\leq t', t\leq 1$ so that
\[
\op{t}{p}=D_1^{m-3}\big(D_1^3-(D_2^2+x_2^2 D_n^2)D_1-{\bar b}\,x_2^3D_n^3\big)=D_1^{m-3}\op{0}{P_{\,\bar b}}.
\]
For any given $u_j(x')\in \gamma^{(s)}$, $j=0, 1, 2$ we define $u_{j+1}(x')\in \gamma^{(s)}$ by setting $u_{j+3}=(D_2^2+x_2^2D_n^2)u_{j+1}+{\bar b}x_2^3D_n^2 u_j$, $j=0,\ldots, m-4$ successively. Assume that  $u\in C^m(\omega)$ satisfies $D_1^{m-3}\op{0}{P_{\,\bar b}}u=0$ with $D_1^ju(0, x')=u_j(x')$, $0\leq j\leq m-1$ then $w=\op{0}{P_{\,\bar b}}$ satisfies $D_1^{m-3}w=0$ with $D_1^j w(0, x')=0$, $j=0,\ldots, m-4$ hence $\op{0}{P_{\,\bar b}}u=0$ contradicting with non local solvability of the Cauchy problem for $\op{0}{P_{\,\bar b}}$.

\medskip
{\it Here is a general question why $s=2$ $($independent of $m\geq 3$$)$}\,?

\medskip
A similar phenomenon is observed in the Cauchy problem for uniformly diagonalizable first order systems (\cite[Theorem 3.3]{Ni:3}). When $m=2$ we have a result similar to Proposition \ref{pro:m:3}:
\begin{prop}
\label{pro:m:2} Let $n\geq 3$ and consider
\[
P_{mod}(x,\xi)=\xi_1^2-2x_2\xi_1\xi_n-\xi_2^2-x_2^3\xi_n^2
\]
which is a homogeneous polynomial in $\xi$ of degree $2$ with polynomial coefficients. For any $0\leq t\leq 1$ the Cauchy problem for $\op{t}{P_{mod}}$ is not locally solvable at the origin in $\gamma^{(s)}$ for $s>5$, in particular, ill-posed in $\gamma^{(s)}$, $s>5$ near the origin.
\end{prop}
This result for $\op{0}{P_{mod}}$ was proved in  \cite{BeNi:1}  (where there is some insufficient part of the proof, see the correction given in \cite{Ni:book}). Then to conclude Proposition \ref{pro:m:2} it is enough to note $\op{t}{P_{mod}}=\op{0}{P_{mod}}$ for $0\leq t\leq 1$. 

\medskip
{\it Now we would ask ourselves is there  an example of a homogeneous polynomial $p$  in $\xi=(\xi_1,\xi_2,\ldots, \xi_n)$, $n\geq 3$  of degree $2$ with real analytic coefficients satisfying \eqref{eq:nonchar} and \eqref{eq:hyperb} such that the Cauchy problem for $\op{t}{p}$,  for any $0\leq t\leq 1$,  is ill-posed in $\gamma^{(s)}$, $s>2$.}

\medskip

For the special case $n=2$ ($m=2$) we have;
\begin{prop}
\label{pro:2:2}Consider 
\[
P(x, \xi)=\xi_1^2-2a(x)\xi_1\xi_2+b(x)\xi_2^2+c(x),\quad x=(x_1, x_2)\in\R^2
\]
which is a polynomial in $\xi=(\xi_1,\xi_2)$ of degree $2$ without homogeneous part of degree $1$ with real analytic coefficients $a(x), b(x), c(x)$ such that $\Delta(x)=a^2(x)-b(x)\geq 0$ near the origin $($$a(x)$, $b(x)$ are assumed to be  real$)$. Then the Cauchy problem for $\op{w}{P}$ is $C^{\infty}$ well-posed near the origin.
\end{prop}
In fact if we make a real analytic change of coordinates $y=\kappa(x)=(x_1, \phi(x))$ such that 
$\phi_{x_1}(x)-a(x)\phi_{x_2}(x)=0$, $\phi(0, x_2)=x_2$ where $\phi_{x_j}=\dif\phi(x)/\dif x_j$ we see that
\begin{equation}
\label{eq:henkan}
\begin{split}
&\op{w}{P(x,\xi)}(u\circ \kappa)\\
=&\Big({\rm op}^0\big(\eta_1^2-\al{\tilde \Delta}\eta_2^2+\be_1{\tilde \Delta_{x_2}}\eta_2+\be_2{\tilde \Delta}\eta_2+\be_3\eta_1+\be_4\big)u\Big)\circ \kappa
\end{split}
\end{equation}
where ${\tilde \Delta}=\Delta\circ\kappa^{-1}$, ${\tilde \Delta_{x_2}}={ \Delta_{x_2}}\circ\kappa^{-1}$ and $\al=\al(y)\geq c>0$, $\be_i=\be_i(y)$ are real analytic near $y=0$. To prove the result, noting that $\big|{\tilde \Delta_{x_2}}\big|\leq C\big|{\tilde \Delta_{y_2}}\big|\leq C'\big|\sqrt{{\tilde \Delta}}\big|$,   it suffices to apply  \cite[Theorem 1.1]{Ni:2} to the right-hand side of \eqref{eq:henkan}.

\medskip

 {\it We could possibly ask if there exists some ${\bar s}>2$ such that for any homogeneous polynomial $p$ in $\xi$ of degree $2$ with real analytic coefficients satisfying \eqref{eq:nonchar} and \eqref{eq:hyperb} the Cauchy problem for $\op{ w}{p}$ is well-posed in $\gamma^{(s)}$ for $s<{\bar s}$ near the origin.}


\end{document}